\documentclass[letterpaper, 10pt, conference]{IEEEtran}

\IEEEoverridecommandlockouts                              

\usepackage[bmargin=25.4mm, rmargin=19.1mm, lmargin=19.1mm, tmargin=19.1mm]{geometry}
\usepackage{graphicx}
\usepackage{subcaption}
\usepackage{siunitx}
\usepackage{amsmath}
\usepackage{amssymb}

\usepackage{dsfont}
\usepackage{tikz}

\usepackage{amsthm}
\newtheorem{theorem}{Theorem}

\newtheorem{lemma}[theorem]{Lemma}
\newtheorem{proposition}[theorem]{Proposition}

\newtheorem{example}{Example}
\newtheorem{remark}{Remark}

\usepackage{xcolor}

\definecolor{cof}{RGB}{219,144,71}
\definecolor{pur}{RGB}{186,146,162}
\definecolor{greeo}{RGB}{91,173,69}
\definecolor{greet}{RGB}{52,111,72}

\newcommand{\bal}[1] {\ensuremath{\left(\begin{array}{#1}}}
\newcommand{\ear} {\ensuremath{\end{array}\right)}}

\newcommand{\bals}[1] {\ensuremath{\left[\begin{array}{#1}}} 
\newcommand{\ears} {\ensuremath{\end{array} \right] }} 

\DeclareMathOperator*{\bigtimes}{\raisebox{-0.3ex}{\text{\Large$\times$}}}

\let\leq\leqslant
\let\geq\geqslant

\newcommand{\PP}{\ensuremath{\mathbf{P}}}
\newcommand{\PQ}{\ensuremath{\mathbf{Q}}}

\newcommand{\calA}{\ensuremath{\mathcal{A}}}
\newcommand{\calB}{\ensuremath{\mathcal{B}}}
\newcommand{\calC}{\ensuremath{\mathcal{C}}}

\newcommand{\calF}{\ensuremath{\mathcal{F}}}

\newcommand{\calI}{\ensuremath{\mathcal{I}}}

\newcommand{\calL}{\ensuremath{\mathcal{L}}}

\newcommand{\calS}{\ensuremath{\mathcal{S}}}

\newcommand{\bmat}{\begin{matrix}}
\newcommand{\emat}{\end{matrix}}
\newcommand{\bbm}{\begin{bmatrix}}
\newcommand{\ebm}{\end{bmatrix}}
\newcommand{\bpm}{\begin{pmatrix}}
\newcommand{\epm}{\end{pmatrix}}
\newcommand{\bse}{\begin{subequations}}
\newcommand{\ese}{\end{subequations}}
\newcommand{\beq}{\begin{equation}}
\newcommand{\eeq}{\end{equation}}
\newcommand{\ben}{\begin{enumerate}}
\newcommand{\een}{\end{enumerate}}
\newcommand{\beni}{\renewcommand{\labelenumi}{\roman{enumi}.}
\renewcommand{\theenumi}{\roman{enumi}}\begin{enumerate}}
\newcommand{\eeni}{\end{enumerate}\renewcommand{\labelenumi}{\arabic{enumi}.}
\renewcommand{\theenumi}{\arabic{enumi}}}
\newcommand{\bena}{\renewcommand{\labelenumi}{\alpha{enumi}.}
\renewcommand{\theenumi}{\alpha{enumi}}\begin{enumerate}}
\newcommand{\eena}{\end{enumerate}\renewcommand{\labelenumi}{\arabic{enumi}.}
\renewcommand{\theenumi}{\arabic{enumi}}}
\newcommand{\bit}{\begin{itemize}}
\newcommand{\eit}{\end{itemize}}

\newcommand{\R}{\ensuremath{\mathbb R}}

\newcommand{\N}{\ensuremath{\mathbb N}}

\newcommand{\E}{\ensuremath{\mathbb E}}

\title{\LARGE \bf
Steady-state Analysis of a Neural-cognition Based Human-social Behavior Model}

\author{Jieqiang Wei, Ehsan Nekouei, Junfeng Wu, Vladimir Cvetkovic, Karl H. Johansson  
\thanks{*This work is supported by Knut and Alice Wallenberg Foundation, Swedish Research Council, and Swedish Foundation for Strategic Research.}
\thanks{Jieqiang Wei, Ehsan Nekouei, Karl H. Johansson are with the Department of Automatic Control, School of Electrical Engineering and computer science. 
 KTH Royal Institute of Technology,
 SE-100 44 Stockholm, Sweden.
 {\tt\small \{jieqiang, nekouei, kallej\}@kth.se}.
 \newline
 Junfeng Wu is with College of Control Science and Engineering, Zhejiang University, Hangzhou, China. 
 {\tt\small jfwu@zju.edu.cn}.
 \newline
 V. Cvetkovic is with School of Architecture and the built environment.
 {\tt\small vdc@kth.se}.
}}

\begin{document}
\maketitle
\thispagestyle{empty}
\pagestyle{empty}

\begin{abstract}\label{s:Abstract}
We consider an extension of the Rescorla-Wagner model which bridges the gap between conditioning and learning on a neural-cognitive, individual psychological level, and the social population level. In this model, the interaction among individuals is captured by a Markov process. The resulting human-social behavior model is a recurrent iterated function systems which  behaves differently  from the classical Rescorla-Wagner model due to randomness. Convergence and ergodicity properties of the internal states of agents in the proposed model are studied. 
%
\end{abstract}


\section{Introduction}\label{s:Introduction}


\textcolor{black}{Internal state of a human agent has significant effects on her decision making process.} This state can be associated with a bias, an irrational or emotional disposition, \textcolor{black}{and} has an important evolutionary role affecting a decision that is presumably based on cognition or a calculation as a rational choice. Recent evidence supports an integrated view of cognition and emotion, the neurological basis \textcolor{black}{which are} high connectivity areas of the brain (hubs) \cite{Pessoa08}. \textcolor{black}{In other words, any decision made by a human agent integrates rational (cognitive) and irrational (emotional) components (or dispositions) on a neurological level \cite{kahneman2000}. Thus, the impact of emotions (or bias) needs to be somehow accounted for in any decision of a human agent.}

Recently, Epstein \cite{epstein2014agent_zero} used the Rescorla-Wagner model, see e.g., \cite{Culver2015,Pearce1980,RW1972,Siegel1996}, to study social behavior. The central concept in the work of \cite{epstein2014agent_zero} is the notion of emotional disposition that is based on conditioning. Epstein models the decisions and actions of  agents as a dynamic process in space and time, where decisions to act are bimodal, i.e., an agents either acts or does not act. The trigger for action is surpassing a specified threshold by the combined emotional and rational dispositions which change in time. The most interesting aspect of Epstein's generic study was to show how important the mutual interactions between agents are for social behavior. Interactions between agents has in fact been the prime focus in studies of opinion formation and consensus in social networks, e.g., \cite{Degroot_74,Xia2011,FriedkinJohnsen1999,friedkin2011social,Friedkin321,Tempo2015} where agents are perceived as essentially rational.


\textcolor{black}{In Epstein's model, the interaction topology between agents is time-invariant. However, social studies suggest that the inter-personal interaction topology among agents in a social network is time-varying and possibly random. For example, the author in \cite{Wasserman1980} proposed a social network model, with continuous-time Markovian interaction networks, which is verified by experimental data. The interested reader is referred to \cite{Robins2001,Courville2006,Sajadi2018} and references within for more details on the time-varying interaction topologies in social networks. Motivated by this observation, the current manuscript proposes a generalization of the Rescorla-Wagner wherein the interaction topology among agents is governed by a Markov chain, namely Markovian random graphs.  In this model, the state of the each agent is updated based on the current state, the states from the neighbors, and the external stimuli. It is shown that this model contains many well-known social network models, e.g., Friedkin-Johnsen model \cite{friedkin2011social} and opinion dynamics \cite{Xia2011}, as special examples.  }

For the proposed model, which is a stochastic process, we then prove the convergence of it. More precisely, we distinguish the convergence for the forward and backward process, respectively.  Furthermore, to study the steady-state of the behavior of the process, an ergodic property is obtained. Comparing to \cite{Barnsley1989}, which is closely related to our model, we extend the ergodicity of the process from bounded functions on the Euclidean space to unbounded ones. In an early study \cite{Wei2017} of the proposed model in this paper, the mean square stability was proved. 

\subsection{Paper outline and notations}

The structure of the paper is as follows. In Section~\ref{s:Preliminaries}, we introduce some terminologies and notations. In Section~\ref{s:ourmodel}, a human-social behavior model is proposed based on the Rescorla-Wagner model. Then the convergence and ergodicity of the proposed model are considered in Section \ref{ss:converge} and Section \ref{ss:ergo}, respectively. Discussion and conclusion are given in Section \ref{s:conclusion}.

\textbf{Notations.} The notations used in this paper is fairly standard. With $\R_-,\R_+, \R_{\geq 0}$ and $\R_{\leqslant 0}$ we denote the sets of negative, positive, non-negative, non-positive  real numbers, respectively. $\mathds{1}_n$ is the $n$-dimensional vector containing only ones. We omit the subscript when there is no confusion. $\delta_{ij}$ denotes the Kronecker delta. 
$\E$ is the expectation. For any set $A$, the product space $\bigtimes_{i\in\N_0}\Omega_i$ with $\Omega_i=A$ is denoted as $A^{\N_0}.$ For any matrix $M$, the induced norm of $M$ is denoted as $\|M\|$. 


\section{Preliminaries}\label{s:Preliminaries}

Given a square matrix $A=(a_{ij})_{i,j=1}^n$, let $\rho(A)$ be its spectral radius. The matrix $A$ is \emph{Schur stable} if $\rho(A)<1$. The matrix is \emph{row stochastic} if $a_{ij}\geq 0$ and $\sum_{j=1}^n a_{ij}=1, \forall i.$ The terminologies about Markov chains are kept consistent with \cite{norris1998markov}.

%

%

\subsection{Markov processes with unique stationary distribution}

The definitions for ergodicity is consistent with \cite{probabilityklenke}.

Let $I=\N_0$, a stochastic process $X=(X_t)_{t\in I}$ is called first-order stationary if 
\begin{align*}
\calL[(X_{t+s})_{t\in I}]=\calL[(X_{t})_{t\in I}], \forall s\in I
\end{align*}
where $\calL[Y]$ is the distribution of the random variable $Y$.

Let $X=(X_n)_{n\in \N_0}$ be a stochastic process with values in a Polish space $E$. Without loss of generality, we assume that $X$ is the canonical process on the probability space $(\Omega, \calA, \mathbf{P}) = (E^{\N_{0}}, \calB(E)^{\otimes \N_0}, \mathbf{P}).$ Define the shift operator
\begin{align*}
\tau:\Omega\rightarrow\Omega, (\omega_n)_{n\in\N_0} \mapsto ((\omega_{n+1})_{n\in\N_0}).
\end{align*}
An event $A\in\calA$ is called invariant if $\tau^{-1}(A)=A$. Denote the $\sigma-$algebra of invariant events by 
\begin{align*}
\calI=\{A\in\calA \mid \tau^{-1}(A)=A \}.
\end{align*}
A $\sigma-$algebra $\calI$ is called $\mathbf{P}-$trivial if $\mathbf{P}[A]\in\{0,1\}$ for every $A\in\calI$. The map $\tau$ is called measure preserving if 
\begin{align*}
\mathbf{P}[\tau^{-1}(A)] = \mathbf{P}[A], \forall A\in\calA.
\end{align*}
In this case, $(\Omega, \calA, \mathbf{P}, \tau)$ is called a measure preserving dynamical system. If $\tau$ is measure preserving and $\calI$ is $\mathbf{P}-$trivial, then $(\Omega, \calA, \mathbf{P}, \tau)$ is called ergodic.
Denote $X_n(\omega) = X_0(\tau^n(\omega))$, where $X_0:\omega\rightarrow E$ is the initial distribution. $X$ is stationary if and only if $(\Omega, \calA, \mathbf{P}, \tau)$ is a measure preserving dynamical system. The stochastic process $X$ is called ergodic if $(\Omega, \calA, \mathbf{P}, \tau)$ is ergodic,where $\tau$ is the shift operator.

\begin{theorem}[Individual ergodic theorem, Birkhoff (1931)]
Suppose that $\tau$ is measure-preserving on $(\Omega,\calA,\PP)$ and that $X_0$ is measurable and integrable. Then 
\begin{align}
\lim_{n\rightarrow\infty} \frac{1}{n} \sum_{k=0}^{n-1} X_k = \E[X_0\mid \calI]
\end{align}
with probability 1.  If $\tau$ is ergodic, $\lim_{n\rightarrow\infty} \frac{1}{n} \sum_{k=0}^{n-1} X_k = \E[X_0]$ with probability 1.
\end{theorem}

\section{Human-social Behavior Model Based on Neural Cognition}\label{s:ourmodel}

\textcolor{black}{In this section, we first explain the Rescorla-Wagner model. Then, an extension of this model is introduced in which the interaction topology between agents is derived by a Markov process.}
\subsection{Rescorla-Wagner Model}
One of the most well-known models in \emph{Pavlovian} theory of reinforcement learning, called \emph{Rescorla-Wagner model}, was proposed in \cite{RW1972}. In the classic Rescorla-Wagner model, the conditional stimulus has an associative value $x\in\R$, supposed to be proportional to the amplitude of the conditional response or to the proportion of conditional response triggered by the conditional stimulus. A typical Pavlovian conditioning session is a succession of several trials. Each trial is composed of the presentation of the conditional stimulus followed by the presentation of the unconditional stimulus. On each trial $k$, the associative value of the conditional stimulus are updated according to the following equation 
\begin{align}\label{e:RW}
x(k+1) = x(k)+ \alpha \Big(r(k)- x(k) \Big), 
\end{align}
where $r(k)\in\R$ is the intensity of the unconditional stimulus on that trial, $\alpha$ is the learning parameter, and $x(k)\in\R$ is the associative strength between the conditional stimulus and the unconditional stimulus. Applications of the Rescorla-Wagner model in machine learning, especially Q-Learning, can be found in e.g., \cite{Hampton8360}, where r(t) is the reward which can be modeled as a Markov chain.

In more general Rescorla-Wagner models, more conditional stimulus can be incorporated. Each conditional stimulus has an associative value $x_i$, which is the associative strength of the $i$th conditional stimulus and the unconditional stimulus, namely some degree to which the conditional stimulus alone elicits the unconditional response. The associative value of all the conditional stimuluses are updated according to the following equation
\begin{align}
x_i(k+1) = x_i(k)+ \alpha^i \Big(r(k) & -\sum_{j=1}^n W_{ij} x_j(k) \Big), \nonumber\\
& i=1,\ldots,n, \label{e:RW_general}
\end{align}
where $n$ is the number of conditional stimulus on that trial, $W_{ij}=\frac{1}{n}$ for all $i$ and $j$, $\alpha^i$ is the learning rate of $i$th conditional stimulus. Rescorla-Wagner model is especially successful in explaining the block phenomenon in Pavlovian conditioning with experimental supports \cite{Siegel1996}.  

The Rescorla-Wagner model has been applied to various levels of human behavior that typically involve emotions and conditioning. A study of human-social behavior that based on Rescorla-Wagner model, which was presented in \cite{epstein2014agent_zero}, establish the connection between neural cognition and human behavior in social networks. 

Consider a society composed by $n$ agents denoted $\calS :=\{1,\ldots,n\}.$ One methodology in \cite{epstein2014agent_zero} of describing human behavior is proposed by separating the human psychology into irrational, rational and social parts. The irrational component evolves according the Rescorla-Wagner model \eqref{e:RW_general} with $W_{ij} = \delta_{ij}$.
Here $x_i(k)$ is the irrational component of $i$th agent state, which can be a belief or an opinion depends on the consider scenario, $r(k)\in \{ 0,1 \}$ is a random binary variable, which takes value one for emotion acquisition, and zero for emotion extinction. Then in the model proposed by \cite{epstein2014agent_zero}, human action depends on whether the summation of irrational and rational components of each agent and these of the related neighbors is larger than a given threshold. \textcolor{black}{In what follows, we refer to  $x_i(\cdot)$ as the internal state of agent $i$.}
\subsection{Rescorla-Wagner Model With Markovian Topology}
 \textcolor{black}{The generalized Rescorla-Wagner model with random time-varying topology is given by}
\begin{align}\label{e:our model}
x(k+1) = B_{i_k} x(k) + A (r(k)-W_{i_k} x(k)),
\end{align}
where $x(k)$ is the vector of the state of the agents, learning rate $A$ is a diagonal matrix satisfying $0\leq A \leq I$, $r(k)$ and $i_k$ are Markov chains with finite states, and for each realization of $i_k$, $B_{i_k}$ and $W_{i_k}$ are row-stochastic matrices. The initial condition is set to be $x(0)=x_0\in\R^n$. Here the matrices $B$ and $W$, corresponding to topologies, can incorporate the time-varying networks. Now we can write our model \eqref{e:our model} into a compact form
\begin{align}\label{e:our model_compact}
x(k) & = F_{i_k}x(k-1) + Ar(k) 
\end{align}
where $F_{i_k}= B_{i_k}- AW_{i_k}$. 

The model \eqref{e:our model} include several established models as special cases, which can be seen by the following examples. First, we establish the resemblance of system \eqref{e:our model} with Rescorla-Wagner model.

\begin{example}
It is straightforward to see that the system \eqref{e:our model} is equivalent to Rescorla-Wagner model \eqref{e:RW_general} by taking  $B_{i_k} = I$ and $W_{i_k} = \frac{1}{n}\mathds{1}\mathds{1}^\top$ for all $k$.
\end{example}

Next, the system \eqref{e:our model} is equivalent to some social network model by specifying appropriate parameters.

\begin{example}
Friedkin-Johnsen model \cite{friedkin2011social} captures the opinion dynamics with heterogeneity, i.e., agents can factor their initial opinions (or prejudices) into every iteration of opinion, as follows
\begin{align}\label{e:FJ}
x(k+1) = \Lambda W x(k)+ (I-\Lambda)u, \quad x(0)=u.
\end{align}
where $W$ is a row stochastic matrix, $\Lambda$ is a diagonal matrix satisfying $0\leq \Lambda \leq I$, and $u$ is the initial opinion. It can bee seen that, by setting $A=I-\Lambda$, $B=W$ and $r(k)$ being deterministic and identical for all $k$, system \eqref{e:our model} includes \eqref{e:FJ} as a special case.
\end{example}

\begin{example}
Agreement and disagreement has been an important topic in the study of social networks, see e.g., \cite{ACEMOGLU2010,DeMarzo,Kar2011,Shi2013} and the references within. As an example, the model considered in \cite{Shi2013} is a special case of the considered model. In fact, there are three events for the iterative update for agent $i$, namely attraction, neglect and repulsion, and each of these events can be formulated into \eqref{e:our model} by choosing appropriate parameters. Furthermore, in \cite{Shi2013}, it is assumed that, at each step, one of these three events is chosen randomly accordingly to a given probability. This is a special case of Markov Chains. 
\end{example}

In the following section, we shall study the convergence of the stochastic process \eqref{e:our model_compact}.


%


\section{Convergence and Ergodicity of Internal States in Markovian Rescorla-Wagner Model}\label{s:main}

In this section, we shall study the convergence and ergodicity property, in Section \ref{ss:converge} and \ref{ss:ergo}, respectively, of the model \eqref{e:our model}. We first introduce the following notations which will be used for the analysis.

Let $(\R^n,d)$ be a complete separable locally compact metric space. Let $i\in\calI : = \{1,2,\ldots,N\}$ and $i_0, i_1, \ldots$ be a Markov chain in $\calI$ with probability transition matrix $P=[p_{ij}]$. Denote the right-hand-side of the system \eqref{e:our model_compact} as $w_{i_k}:\R^n\rightarrow\R^n$, i.e., $w_{i_k}(x) = F_{i_k}x+Ar(k)$, which is Lipschitz continuous with respect to $x$.  Assume that the Markov chain $i_0, i_1, \ldots$ admits a unique stationary distribution. Consider a random walk given as 
\begin{equation}\label{e:our model_general}
\begin{aligned}
Z_k & = w_{i_k} (Z_{k-1}) \\
& = w_{i_k}\cdots w_{i_1} (Z_0) \\
& = F_{i_k}\cdots F_{i_1} Z_0 + Ar(k) \\
& \quad + \sum_{\ell = 2}^{k} F_{i_k}\cdots F_{i_\ell}A r(\ell-1)
\end{aligned}
\end{equation}
where $w_{i_k}\cdots w_{i_1}$ denotes the composition of the function sequence, $Z_0:\Omega_c\rightarrow\R^n$ is a random variable on probability space $(\Omega_c,\calF_c,\PP_c)$ defined on $\R^n$.
In this section, we focus on the convergence and ergodicity of the family of random variables $Z = (Z_t, t\in\N_0)$. 

The Markov chain $\{i_k\}$ is defined on the probability space $(\Omega_d,\calF_d,\PP_d)$. The distribution (probability measure) of $i_k$ and $Z_k$ are $\PP_{i_k} = \PP_d\circ i_k^{-1}$ and $\PP_{Z_k} = \PP_c\circ Z_k^{-1}$, respectively. Now the family of random variables $Z = (Z_t, t\in\N_0)$ is a stochastic process on $((\Omega_d\times \Omega_c)^{\N_0}, (\calF_d\times \calF_c)^{\otimes \N_0})$ with value in $(\R^n, \calB(\R^n))$.

\subsection{Convergence of forward and backward processes}\label{ss:converge}

In this subsection, we present some convergence results for the process \eqref{e:our model_general}. Notice that the process $\{Z_k\}$ given by \eqref{e:our model_general} is not a Markov process, but $\tilde{Z}_k:=(Z_k, i_k)$ is a Markov process with values in $\tilde{X}:= \R^n\times \calI$. 

In order to derive the convergence of the distribution of $\{Z_k\}$, one essential part of the techniques is related to the associated backward process $\overleftarrow{Z}$, i.e.,
\begin{align}\label{e:backward}
\overleftarrow{Z}_k  = & w_{i_1}\cdots w_{i_k}Z_0 \\
= & F_{i_1}\cdots F_{i_k} Z_0 + Ar(1) \\
& + \sum_{\ell = 1}^{k-1} F_{i_1}\cdots F_{i_\ell}A r(\ell+1).
\end{align}

Let $(m_i)$ be the unique stationary initial distribution for the Markov chain $i_0, i_1, \ldots$ on $\calI$, i.e.,
\begin{align}
\sum_{i=1}^N m_i p_{ij} = m_j, \quad j= 1,\ldots,N.
\end{align}
Let the matrix $Q=[q_{ij}] \in\R^{N\times N}$, given as 
\begin{align}
q_{ij} = \frac{m_j}{m_i}p_{ji},
\end{align}
be the inverse transition probability matrix (see e.g., Theorem 1.9.1 in  \cite{norris1998markov}), namely the probability that $(i_1,\ldots,i_k) = (j_1,\ldots,j_k)$ is 
\begin{align}
\sum_{j(0)=1}^N m_{j_0}p_{j_0j_1}\cdots p_{j_{k-1}j_k} 
=  m_{j_k}q_{j_kj_{k-1}}\cdots q_{j_2j_1}.
\end{align}
Denote the set $\Omega = \{\mathbf{i} = (i_0,i_1, \ldots) \} = \calI^{\N_0}.$ Let $\PP$ be the probability on $\Omega$ for the forward chain $\PP(i_0,i_1, \ldots, i_k) = m_{i_0}p_{i_0i_1}\cdots p_{i_{k-1}i_k}$ and $\PQ$ be the probability corresponding to the backward chain $\PQ(i_0,i_1, \ldots, i_k) = m_{i_0}q_{i_0i_1}\cdots q_{i_{k-1}i_k}$.

We first recall a result, presented in \cite{Barnsley1989}, which gives the convergence of the forward and backward process, respectively. Notice that the behaviors of forward and backward processes are very different in the sense that forward process converge in distribution while the backward process converges almost surely. Moreover, the convergence result for forward process is for the initial distribution $\tilde{\nu}$ satisfying $\tilde{\nu} (\R^n\times \{i\})= \PP_{i_0} (\{i\}) =m_i, i=1,\ldots, N$, i.e., the Markov chain $\{i_k\}$ is initialized with stationary distribution, while $Z_0$ is arbitrary. Recall that  
\begin{align}
& \E_\PP (\log\| F_{i_k}\cdots F_{i_1} \|) \nonumber\\
= & \sum_{i_1}\cdots \sum_{i_k} m_{i_1}p_{i_1i_2}\cdots p_{i_{k-1}i_k} \log\| F_{i_k}\cdots F_{i_1} \|.
\end{align}

\begin{lemma}[Theorem 2.1, \cite{Barnsley1989}]\label{lm:base}
If, for some $k$, 
\begin{align}\label{e:stablizing-condition}
\E_\PP (\log\| F_{i_k}\cdots F_{i_1} \|) <0,
\end{align}
then 
\begin{enumerate}
\item[(1)] for $\PQ$ almost all $\mathbf{i}$, the backward process $\overleftarrow{x}_k = w_{i_1}\cdots w_{i_k}x_0 $ converges to a random variable, denoted as $Y(\mathbf{i})$, as $k\rightarrow\infty$, which does not depend on $x_0$. In other words, for given $x_0$, the random variable $\overleftarrow{x}_k:\Omega\rightarrow \R^n$ converges to a finite limit $\PQ-$almost surely.
\end{enumerate}
Define the distribution of $(Y, i_1)$ on $\tilde{X}$ as $\tilde{\mu}(\tilde{B}) = \PQ(\mathbf{i}:(Y(\mathbf{i}),i_1(\mathbf{i}))\in\tilde{B})$, where $\tilde{B}\subset \tilde{X}$ is a Borel set,
\begin{enumerate}
\item[(2)] then $\tilde{\mu}$ is the unique stationary initial distribution for the Markov process $(Z_k,i_k)$; \\
\item[(3)] for any probability measure on $\tilde{X}$, denoted as $\tilde{\nu}$, satisfying $\tilde{\nu} (\R^n\times \{i\})=m_i, i=1,\ldots, N$, then the random walk $\tilde{Z}_k^{\tilde{\nu}}$, i.e., the Markov process with initial distribution $\tilde{\nu}$, converges in distribution to $\tilde{\mu}$. Furthermore, the random walk $Z_k^{\tilde{\nu}}$ on $\R^n$ converges in distribution to the measure $\mu(B) = \tilde{\mu}(B\times \calI)$. 
\end{enumerate}
\end{lemma}

The previous lemma shows that if the initial distribution corresponding to $i_0$ is stationary, then $Z_k^{\tilde{\nu}}$ converges to $\mu$ in distribution. In the following result, we extend the  result to arbitrary initial distribution $\tilde{\nu}^0$ for both $Z_0$ and $i_0$. Here $\tilde{\nu}^0(\R^n\times\{i\}) = \eta^0_i$. Denote the distribution of the Markov chain $i_k$ with initial distribution $\eta^0_i$ as $\eta^k_i, i=1,\ldots,N$. Based on the distribution $\eta^k$ at time $k$, we define a probability $\PP'_k$ on $\Omega$ as $\PP'_k(i_0,\ldots,i_n) = \eta^k_{i_0}p_{i_0 i_1}\cdots p_{i_{n-1} i_n}$. The process with initial distribution $\tilde{\nu}^0$ is denoted as $\tilde{Z}^{\tilde{\nu}^0}_k$. Here $\tilde{Z}^{\tilde{\nu}^0}_k$ is a random variable $\R^n\times \Omega\rightarrow \R^n\times \calI$.

\begin{proposition}\label{prop:arbit-initial}
Assume that $\underline{m} := \min_i m_i >0$. Then under the same assumptions as in Lemma \ref{lm:base} and $p_{ij}<1$ for $\forall i,j=1,\ldots,N$, the random walk $\tilde{Z}^{\tilde{\nu}^0}_k$ converges in distribution to $\tilde{\mu}$.
\end{proposition}
\textcolor{black}{Proposition \ref{prop:arbit-initial} implies that the distribution of the internal states converges to the stationary distribution, induced by the generalized model, regardless of the initial distribution of the Markov chain generating the interaction topologies.}
\begin{proof}
Let $\tilde{\nu}^0$ be the initial distribution and $\tilde{\nu}^k$ be the distribution of $\tilde{Z}^{\tilde{\nu}^0}_k$.  Since $i_k$ has the unique stationary distribution, we have that $\eta^k_i\rightarrow m_i$ as $k\rightarrow\infty$. Hence, for any $\varepsilon>0$, there exists $M$ such that $\|\eta^k-m\|<\varepsilon$ for any $k>M$. For each $j$ and $k$, let the conditional distribution of $Z_k$ given $i_k=j$ to be denoted as
\begin{align}
\nu^k_j = \frac{\tilde{\nu}^k(B\times \{j\})}{\tilde{\nu}^k(X\times \{j\})}.
\end{align}
For all $\tilde{f}\in \calC_b(\tilde{X})$, continuous and bounded function $\tilde{X}\rightarrow \R$, we have 
\begin{align}
& \E[\tilde{f}(\tilde{Z}^{\tilde{\nu}^0}_k)] \nonumber\\ 
= & \int_{\mathbf{i}=\{i_0, \ldots, i_n \}}\int \tilde{f}(w_{i_n}\cdots w_{i_1}x,i_n) d\nu^0_{i_0}(x) d \PP'_0(\mathbf{i}) \nonumber \\
= & \int_{\mathbf{i}=\{i_1, \ldots, i_n \}}\int \tilde{f}(w_{i_n}\cdots w_{i_2}x,i_n) d\nu^1_{i_1}(x) d \PP'_1(\mathbf{i}) \nonumber \\ 
= & \int_{\mathbf{i}=\{i_M, \ldots, i_n \}}\int \tilde{f}(w_{i_n}\cdots w_{i_{M+1}}x,i_n) d\nu^M_{i_M}(x) d \PP'_M(\mathbf{i}) \nonumber \\
\in & \int_{\mathbf{i}=\{i_M, \ldots, i_n \}}\int \tilde{f}(w_{i_n}\cdots w_{i_{M+1}}x,i_n) d\nu^M_{i_M}(x) d \PP(\mathbf{i})  \nonumber \\
& + B(0,\varepsilon\tilde{K}) \nonumber
\end{align}
where the last inclusion is based on the following derivations. First, notice that  
\begin{align*}
& \int_{\mathbf{i}=\{i_M, \ldots, i_n \}}\int \tilde{f}(w_{i_n}\cdots w_{i_{M+1}}x,i_n) d\nu^M_{i_M}(x) d \PP'_M(\mathbf{i}) \\
= & \sum_{\mathbf{i}=\{i_{M+1}, \ldots, i_n \}} \big( \sum_{i_M}\int \tilde{f}(w_{i_n}\cdots w_{i_{M+1}}x,i_n) d\nu^M_{i_M} \\   
& \quad \quad \quad \quad \eta^M_{i_M} p_{i_M i_{M+1}}  \big) p_{i_{M+1} i_{M+2}}\cdots p_{i_{n-1} i_{n}}  \\
\in & \sum_{\mathbf{i}=\{i_{M+1}, \ldots, i_n \}} \big( \sum_{i_M}\int \tilde{f}(w_{i_n}\cdots w_{i_{M+1}}x,i_n) d\nu^M_{i_M}   \\
& \quad \quad \quad \quad (m_{i_M}+B(0,\varepsilon)) p_{i_M i_{M+1}}  \big) p_{i_{M+1} i_{M+2}}\cdots p_{i_{n-1} i_{n}}.
\end{align*}
Furthermore, since 
\begin{align*}
& \sum_{\mathbf{i}=\{i_{M+1}, \ldots, i_n \}} \big( \sum_{i_M}\int \tilde{f}(w_{i_n}\cdots w_{i_{M+1}}x,i_n) d\nu^M_{i_M}  \\ 
& \quad \quad \quad \quad \quad m_{i_M} p_{i_M i_{M+1}}  \big)  p_{i_{M+1} i_{M+2}}\cdots p_{i_{n-1}} \\
= & \int_{\mathbf{i}=\{i_M, \ldots, i_n \}}\int \tilde{f}(w_{i_n}\cdots w_{i_{M+1}}x,i_n) d\nu^M_{i_M}(x) d \PP(\mathbf{i})
\end{align*}
and
\begin{align}
& \sum_{\mathbf{i}=\{i_{M+1}, \ldots, i_n \}} \big( \varepsilon \sum_{i_M}\int \tilde{f}(w_{i_n}\cdots w_{i_{M+1}}x,i_n) d\nu^M_{i_M}   \nonumber\\
& \quad \quad \quad \quad \quad  p_{i_M i_{M+1}}  \big) p_{i_{M+1} i_{M+2}}\cdots p_{i_{n-1} i_{n}} \nonumber\\
\leq & \sum_{\mathbf{i}=\{i_{M+1}, \ldots, i_n \}} (K\varepsilon N \bar{P})p_{i_{M+1} i_{M+2}}\cdots p_{i_{n-1} i_{n}} \label{e:ineq-PP}
\end{align}
where $\bar{P} = \max_{i,j\in \calI} p_{ij}$ and $|f|<K$. Moreover, since $\sum_{i_{M+1},\cdots, i_n} m_{M+1} p_{i_{M+1}i_{M+2}}\cdots p_{i_{n-1}i_n} =1$, we have \eqref{e:ineq-PP} is no bigger than $\frac{K\varepsilon N \bar{P}}{\underline{m}}:= \varepsilon\tilde{K}$ where $\underline{m} = \min_i m_i >0$. 

The rest of the proof is based on \cite{Barnsley1989}. Since
\begin{align*}
& \int_{\mathbf{i}=\{i_M, \ldots, i_n \}}\int \tilde{f}(w_{i_n}\cdots w_{i_{M+1}}x,i_n) d\nu^M_{i_M}(x) d \PP(\mathbf{i}) \\
= & \int_{\mathbf{i}=\{i_M, \ldots, i_n \}}\int \tilde{f}(w_{i_{M+1}}\cdots w_{i_n}x,i_n) d\nu^{n+1}_{i_{n+1}}(x) d \PQ(\mathbf{i}) \nonumber 
\end{align*}
which converges to 
\begin{align*}
& \lim_{n\rightarrow\infty} \int_{\mathbf{i}=\{i_M, \ldots, i_n \}}\int \tilde{f}(w_{i_{M+1}}\cdots w_{i_n}x_0,i_n) \\ 
& \quad \quad \quad \quad \quad \quad \quad \quad  d\nu^{n+1}_{i_{n+1}}(x) d \PQ(\mathbf{i}) \\
= & \lim_{n\rightarrow\infty} \int_{\mathbf{i}=\{i_M, \ldots, i_n \}} \tilde{f}(w_{i_{M+1}}\cdots w_{i_n}x_0,i_n) d \PQ(\mathbf{i}) \\
= & \int \tilde{f} \tilde{\mu}, \forall x_0
\end{align*}
then the conclusion follows from Portemanteau's Theorem \cite{probabilityklenke}.
\end{proof}

\subsection{Ergodicity}\label{ss:ergo}

In this section, we present ergodic result about system \eqref{e:our model} which is a version of the strong law of large numbers. 

\begin{theorem}\label{Theo: ERG}
Consider the stochastic process \eqref{e:our model} initialized with arbitrary distribution $\nu$ and the Markov chain $\{i_k\}$ initialized with stationary distribution $(m_i)$.
If, for some $k$, 
\begin{align}
\E_\PP (\log\| F_{i_k}\cdots F_{i_1} \|) <0,
\end{align} 
then 
\begin{align}
\lim_{n\rightarrow\infty} \frac{1}{n} \sum_{k=0}^{n} Z^{\nu}_k = \E[Z^{\mu}_0]
\end{align}
almost surely, where $\mu$ is given in Lemma \ref{lm:base}.
\end{theorem} 
\textcolor{black}{Theorem \ref{Theo: ERG} establishes the ergodic property of the proposed modeled. According to this result, the limiting behavior of the time-average of the internal states in a social network with Markovian interaction topologies can be characterized by the stationary distribution imposed by the dynamics. }
\begin{proof}

Denote the initial distribution of the augmented state as $\tilde{\nu}$. Denote $(\varOmega, \calA) = (\tilde{X}^{\N_{0}}, \calB(\tilde{X})^{\otimes \N_0 })$. Then $\tilde{Z}^{\tilde{\nu}}=(\tilde{Z}^{\tilde{\nu}}_k)_{k\in\N_0}$ is a Markov process with value in $\tilde{X}$. 
Define the shift operator 
\begin{align}
\tau : \varOmega\rightarrow\varOmega, \quad (\omega_n)_{n\in \N_0} \rightarrow (\omega_{n+1})_{n\in \N_0}.
\end{align}
Then $\tilde{Z}^{\tilde{\nu}}_k(\omega) =\tilde{Z}^{\tilde{\nu}}_0(\tau^k(\omega))$. 

First, by Lemma \ref{lm:base} and Corollary 12 in \cite{Hairer2006}, we have that the operator $\tau$ is ergodic. Then by Birkhoff’s Ergodic Theorem \cite{probabilityklenke}, we have 
\begin{align*}
\lim_{n\rightarrow\infty} \frac{1}{n} \sum_{k=0}^{n} \tilde{Z}^{\tilde{\mu}}_k = \E[\tilde{Z}^{\tilde{\mu}}_0]
\end{align*}
which implies that 
\begin{align*}
\lim_{n\rightarrow\infty} \frac{1}{n} \sum_{k=0}^{n} Z^{\mu}_k = \E[Z^{\mu}_0].
\end{align*}
Notice that
\begin{align}
& \|\frac{1}{n}\sum_{k = 0}^{n-1} Z^{\nu}_k -  \E[Z^{\mu}_0] \|_1 \\
= & \| \frac{1}{n}\sum_{k = 0}^{n-1} (Z^{\nu}_k - Z^{\mu}_k ) + \frac{1}{n}\sum_{k = 0}^{n-1} Z^{\mu}_k -  \E[Z^{\mu}_0] \|_1 \\
\leq & \|\frac{1}{n}\sum_{k = 0}^{n-1} (Z^{\nu}_k - Z^{\mu}_k ) \|_1 + \|\frac{1}{n}\sum_{k = 0}^{n-1} Z^{\mu}_k -  \E[Z^{\mu}_0] \|_1.
\end{align}
Moreover, 
\begin{align}
& \|\frac{1}{n}\sum_{k = 0}^{n-1} (Z^{\nu}_k - Z^{\mu}_k ) \|_1 \\
\leq & \frac{1}{n}\sum_{k = 0}^{n-1} \| Z^{\nu}_k - Z^{\mu}_k \|_1
\end{align}
and for any $\varepsilon>0$
\begin{align}
& \PP_c \Big( \|Z^{\nu}_k - Z^{\mu}_k\|_1 \geq \varepsilon^ k \Big) \\
\leq & \frac{\E(\|Z^{\nu}_k - Z^{\mu}_k\|_1)}{\varepsilon^k} \\
\leq & \frac{\E(\|F_k \cdots F_1\| \| (Z^{\nu}_0 - Z^{\mu}_0) \|_1)}{\varepsilon^k} \\
\leq & \frac{\E(\|F_k \cdots F_1\|) \E ( \| (Z^{\nu}_0 - Z^{\mu}_0) \|_1)}{\varepsilon^k}
\end{align}
where the last inequality is implied by H\"{o}lder inequality. Since the Markov chain $\{i_k\}$ initialized with stationary distribution $(m_i)$, it is showed in \cite{Barnsley1989} that condition \eqref{e:stablizing-condition} is equivalent to 
\begin{align}
\lim_{n\rightarrow \infty} \frac{1}{n} \log \|F_{i_n}\cdots F_{i_1}\| = -\alpha, \PP-\textnormal{almost surely},
\end{align}
for $\PP$ almost all $\mathbf{i}$. Hence  $n_0$ can be chosen such that for any $k\geq n_0$ we have $\|w_k \cdots w_1\| < e^{-\frac{k\alpha}{2}}$. Then
\begin{align}
& \PP_c \Big( \|Z^{\nu}_k - Z^{\mu}_k\|_1 \geq \varepsilon^ k \Big) \\
\leq & \frac{e^{-\frac{k\alpha}{2}}}{\varepsilon^ k} \E ( \| (Z^{\nu}_0 - Z^{\mu}_0) \|_1). 
\end{align}
Then if $\varepsilon\in (e^{-\frac{\alpha}{2}},1)$, then the Borel-Cantelli Lemma implies that with probability one $\|Z^{\nu}_k - Z^{\mu}_k\|_1 \leq \varepsilon^ k$ for all but finitely many values of $k$. Therefore, almost surely $ \frac{1}{n} \sum_{k=0}^{n-1} \|Z^{\nu}_k - Z^{\mu}_k\|_1 $ converges to zero as $n\rightarrow\infty$. Hence the conclusion follows.
\end{proof}

%

\begin{remark}
Compared to the result in Theorem 2.1 (iii) \cite{Barnsley1989}, where the ergodicity is proved for the process with bounded continuous function, i.e., $\lim_{n\rightarrow\infty} \frac{1}{n}\sum_{k=0}^{n} f(Z^{\nu}_{k})$ with bounded continuous $f$,  we extend the ergodicity property for identity function which is not bounded.
\end{remark}

\begin{remark}
For the system \eqref{e:our model}, one sufficient condition which guarantees $\E_\PP (\log\| F_{i_k}\cdots F_{i_1} \|) <0$ for some $k$ is that $F_{i}$ is Schur stable for any $i\in\calI$. Notice that the results in this section do not guarantee any boundedness of the states of system \eqref{e:our model_general}. In fact, there are examples satisfying $A_{i}$ is Schur stable for any $i\in\calI$, but the states diverge to infinity with positive probability, see Example 3.17 in \cite{Costa2005}.   
\end{remark}

%

%
%
%
%

\section{Conclusion}\label{s:conclusion}

In this paper, we propose a human-social behavior model, which is based on the well-known Rescorla-Wagner model from neural-cognition and Markovian social networks. The proposed model contains the classical Rescorla-Wagner model and Friedkin-Johnsen model as special cases. Under a sufficient condition, different convergence behaviors for the forward process and backward process are discussed. For the steady-state behavior of the forward process, the ergodicity is proved under the same sufficient condition. 
Incorporation of the proposed model into human decision-making process within a social network is the direction of our future study.

\section{Acknowledgment}

The authors would like to acknowledge Dr. Anton V. Proskurnikov for the constructive discussions.




\bibliographystyle{plain} 
\bibliography{ref}

\end{document}